\documentclass[final,hidelinks,onefignum,onetabnum]{siamart220329}

\usepackage{lipsum}
\usepackage{amsfonts}
\usepackage{graphicx}
\usepackage{epstopdf}
\usepackage{algorithmic}
\usepackage{tikz}
\usepackage{pgfplots}
\usepackage[caption=false]{subfig}
\ifpdf
  \DeclareGraphicsExtensions{.eps,.pdf,.png,.jpg}
\else
  \DeclareGraphicsExtensions{.eps}
\fi

\newsiamremark{remark}{Remark}
\newsiamremark{hypothesis}{Hypothesis}
\crefname{hypothesis}{Hypothesis}{Hypotheses}
\newsiamthm{claim}{Claim}

\headers{Optimisation of Power Grid Stability Under Uncertainty}{John M. Moloney, Samuel J. Williamson, and Cameron L. Hall}

\title{Optimisation of Power Grid Stability Under Uncertainty\thanks{Submitted to the editors DATE.
\funding{This work was funded by the Engineering and Physical Sciences Research Council under grant no.~EP/T517823/1.}}}

\author{John M. Moloney\thanks{Department of Engineering Mathematics, University of Bristol, Bristol, UK 
  (\email{john.moloney@bristol.ac.uk} \email{cameron.hall@bristol.ac.uk}).}
\and Samuel J. Williamson\thanks{Department of Electrical and Electronic Engineering, University of Bristol, Bristol, UK 
  (\email{sam.williamson@bristol.ac.uk}).}
\and Cameron L. Hall\footnotemark[2]}

\usepackage{amsopn}

\ifpdf
\hypersetup{
  pdftitle={Optimisation of Power Grid Stability Under Uncertainty},
  pdfauthor={John M. Moloney, Sam J. Williamson, and Cameron L. Hall}
}
\fi

\begin{document}

\maketitle

\begin{abstract}
The increased integration of intermittent and decentralised forms of power production has eroded the stability margins of power grids and made it more challenging to ensure reliable and secure power transmission. Reliable grid operation requires system-scale stability in response to perturbations in supply or load; previous studies have shown that this can be achieved by tuning the effective damping parameters of the generators in the grid. In this paper, we present and analyse the problem of tuning damping parameters when there is some uncertainty in the underlying system. We show that sophisticated methods that assume no uncertainty can yield results that are less robust than those produced by simpler methods. We define a quantile-based metric of stability that ensures that power grids remain stable even as worst-case scenarios are approached, and we develop optimisation methods for tuning damping parameters to achieve this stability. By comparing optimisation methods that rely on different assumptions, we suggest efficient heuristics for finding parameters that achieve highly stable and robust grids.
\end{abstract}



\section{Introduction}
\label{Introduction}
The urgent need to decarbonise power generation necessitates increasing the amount of power generated from renewable sources. However, this poses new challenges in ensuring the stable and reliable operation of power grids \cite{Shair2021}. As we replace traditional synchronous generators with inverter-based technologies, the support provided to the grid in the form of stored rotational energy is reduced. This lack of stored energy through mechanical inertia inhibits the ability of the grid to damp fluctuations in frequency \cite{Milano2018}. These fluctuations are a result of differences in power supply and demand which are becoming more prevalent as we integrate intermittent renewable power sources into the grid \cite{Ratnam2020}. Fluctuations in frequency cause issues with grid components and may potentially lead to failures or blackouts \cite{Chow2014,Zhang2016}. This paper aims to develop methods for finding the power grid parameters that ensure the grid is inherently stable to frequency fluctuations, even as grids incorporate intermittent and decentralised sources of power into the grid.

One way in which oscillations in frequency are addressed is by providing sufficient damping torque through fast-acting control loops, commonly referred to as power system stabilisers \cite{Devarapalli2022,Kundur1994}. Mathematically this can be treated as optimising the necessary damping required from the power system stabilisers at the generation units. Molnar et al. \cite{Molnar2020,Molnar2021} and Motter et al. \cite{Motter2013} use dynamical systems theory to show how the damping associated with each generator impacts the stability of the system. In particular, Molnar et al. \cite{Molnar2021} identify the effective damping parameter as a key factor in stabilising power grids. This effective damping parameter is composed of a damping parameter that is tuneable using power system stabilisers combined with an inertial parameter that is not tuneable. Research into synthetic and virtual inertia suggests the inertial damping parameter may also be controllable \cite{Tamrakar2017,Kerdphol2019} but this has yet to be implemented in large-scale power grids. 

The methods developed by Molnar et al. \cite{Molnar2021} involve using an optimisation algorithm to determine the necessary electro-mechanical damping parameter for each generator in the power grid. Molnar et al. \cite{Molnar2021} show that by tuning the effective damping parameters appropriately, it is possible to strengthen the stability of the power grid. Although one judiciously-chosen effective damping parameter applied to each generator in the network gives good stability, Molnar et al. \cite{Molnar2021} show that this can be improved further by relaxing the constraint that all effective damping parameters must be the same, instead allowing optimisation based on heterogeneous effective damping parameters.

The results in \cite{Molnar2021} also show that the stability landscape is complex with steep ridges and cusps common throughout the hypersurface. This makes it particularly challenging to find global minimisers and requires sophisticated approaches to optimisation in \cite{Molnar2021}. However, the complexity of the stability landscape also raises the question of robustness; even if we find the global optimum, this solution may rapidly become much worse in the immediate vicinity. If a perturbation is applied to model parameters then the stability may be worse than expected. Given that there is uncertainty in power system parameters, we would like to search the stability landscape and find a solution which is optimal subject to a given level of uncertainty.

Previous work on incorporating uncertainty into power systems analysis has revolved around either the design of power system stabilisers which incorporate uncertainty into power system operation or the choice of an appropriate control algorithm \cite{Chen1995, ElMetwally1995, Bevrani2011}. The control algorithms include robust and adaptive control, both of which are used to damp oscillations while making use of either centralised or distributed power system stabilisers \cite{Lu2015, Tu2019}. Although this effectively addresses the issue of uncertainty in power systems by measuring the state of the system before deciding on appropriate action, this is a reactionary process and does not provide any inherent stability to the power grid. In the event these control systems experience issues, this would leave the grid without any stabilising control. On very short timescales, before control algorithms can be effective, it is necessary to have the power system in as stable a state as possible through tuned effective damping parameters.

It is difficult to choose the optimal effective damping parameters due to the inherent mechanical noise in power system stabilisers which may result in some deviation from the desired damping \cite{Heylen2022, Kerdphol2019}. There is also uncertainty with each measurement taken in a power system due to the tolerances associated with any measurement device \cite{1995a}. These two sources of noise contribute to uncertainty in important power system parameters such as power generation, power demand, conductance, susceptance, inertia and damping. Therefore, there is uncertainty in the power system parameters that determine the steady state power flow and there is uncertainty in the tuneable damping parameters in the sense that they can never be specified precisely, only up to a precision determined by mechanical noise. In this paper, we investigate the influence of these uncertainties associated with a power system on the optimal stability of a power grid.

We first present a mathematical model of power grid stability in \cref{Mathematical Model of Stability} while outlining the parameters of the model and their associated uncertainty. We then describe optimisation methods in \cref{Optimisation Under Uncertainty} that find the effective damping parameters that ensure the best stability of the power grid, provided there is limited noise. We show that as the noise applied to the effective damping increases, the resulting Lyapunov exponent distributions are no longer optimal. We introduce a method which incorporates uncertainty in the optimisation process by sampling from a distribution of effective damping parameters multiple times and taking the average Lyapunov exponent. We run a simulated annealing algorithm using this measure of the Lyapunov exponent under uncertainty and use it to find optimal stability results for this model. We examine a variety of cases where the uncertainty associated with the effective damping parameters is known with various degrees of certainty. In \cref{Statistical Analysis} we then analyse each model statistically using a quantile versus Lyapunov exponent plot to determine the performance of each method when noise is applied to the full power system and to establish which method is appropriate for different levels of uncertainty. Finally, we discuss our results in \cref{Discussion and Conclusions} as well as the implications of this analysis and areas of future work.

\section{Mathematical Model of Stability}
\label{Mathematical Model of Stability}
Molnar et al. \cite{Molnar2021} developed their model of power grid stability based on the swing equation in the following form,
\begin{equation}
	\label{eq:1}
	\ddot{\delta}_{i} + \beta_{i}\dot{\delta}_{i} = \alpha_{i} - \sum_{k \neq i}^{}c_{ik}\sin(\delta_{i} - \delta_{k} - \gamma_{ik}),
\end{equation}
where the phase angle for generator $i$ relative to a reference frequency is given by $\delta_{i}$; the effective damping parameter is given by $\beta_{i} = \frac{D_{i}}{2H_{i}}$ (with $D_{i}$ representing a combination of mechanical and electrical effects that result in damping, and $H_{i}$ representing the inertia constant); the net power driving the generator is given by $\alpha_{i}$; and the coupling strength and the phase shift characterizing the electrical interactions are given by $c_{ik}$ and $\gamma_{ik}$ respectively.

The stability of the synchronous state is then analysed against small perturbations using the Jacobian matrix,
\begin{equation}
	\label{eq:2}
	\mathbf{J} = \begin{pmatrix}
		\mathbf{O} & \mathbf{I}\\
		-\mathbf{P} & -\mathbf{B}
	\end{pmatrix},
\end{equation}
where $\mathbf{O}$ and $\mathbf{I}$ denote the $n \times n$ null and identity matrices respectively. The $n \times n$ matrix $\mathbf{P}$ is defined by
\begin{equation}
	\label{eq:3}
	P_{ik} = 
	\begin{cases}
		-c_{ik}\cos(\delta_{i}^{*} - \delta_{k}^{*} - \gamma_{ik}), \hspace{4.75mm} i \neq k,\\
		\hspace{10mm} -\sum_{k^{'}\neq i}P_{ik^{'}}, \hspace{10mm} i = k,
	\end{cases}
\end{equation}
which expresses the interactions between the generators, while the $n \times n$ diagonal matrix $\mathbf{B}$ has $\beta_{i}$ as its diagonal components. The parameter $\delta_{i}^{*}$ is the steady state phase angle for generator $i$ when synchronised to the nominal frequency. The goal is to choose the effective damping parameters $\{\beta_{i}\}$ to achieve the quickest return to equilibrium following a perturbation. Mathematically, this corresponds to minimising the Lyapunov exponent,
\begin{equation}
	\label{eq:4}
	\lambda_{L} = \max_{i\geq2}\Re(\lambda_{i}) 	
\end{equation}
where $\lambda_{i}$ are the eigenvalues of $\mathbf{J}$. A more negative Lyapunov exponent corresponds to a more stable network. 

One way of determining the optimal $\{\beta_{i}\}$ is to run an optimisation algorithm that searches the solution space of $\lambda_{L}$ for a global minimum; the set of $\beta_{i}$ values found in this way is denoted as $\beta_{\neq}$. A simpler optimisation can be obtained by requiring the $\beta_i$ values to be equal; the single $\beta$ value that minimises $\lambda_{L}$ In this case is denoted $\beta_{=}$. The main advantage of using $\beta_{=}$ is that we reduce the solution space to one dimension, reducing the computational expense of optimisation.

\begin{figure}[h!]
	\centering
	\includegraphics[width=1.0\textwidth]{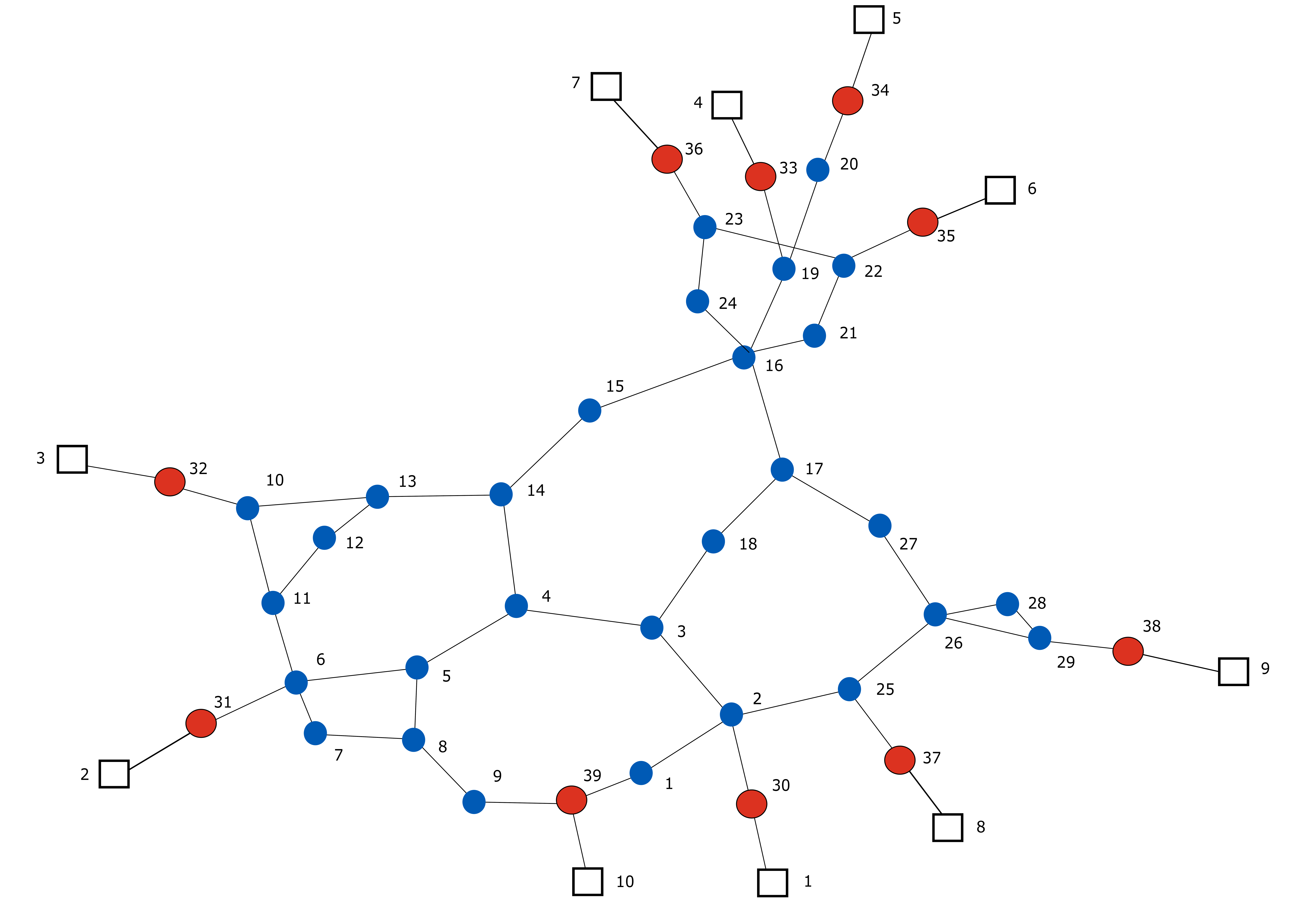}
	\caption{Power grid network used in this analysis where each bus is represented by a circle with the red circles being the generator buses, the blue circles being the load buses and the white squares representing the generators which are connected to a single generator bus.}
	\label{fig:Power_Network}
\end{figure}

In this paper, we use the commonly-studied New England test system to investigate all methods considered. The IEEE 10-generator, 39-node test system is described in \cite{Pai1989} and \cite{Athay1979} as shown in \cref{fig:Power_Network}. We use the data file (case39.m) which can be found in the MATPOWER toolbox \cite{Zimmerman2011} and add dynamic parameters from \cite{Pai1989} similar to those given in \cite{Molnar2021}.

There is uncertainty denoted as $\Delta$ associated with several of the parameters in \cref{eq:1}. The quantities $\Delta \gamma_{ik}$ and $\Delta c_{ik}$ both depend on the uncertainty in the admittance $\Delta Y_{i}$ \cite{Nishikawa2015}. This uncertainty $\Delta Y_{i}$ is a result of uncertainty in the conductance $\Delta G_{i}$ and susceptance $\Delta B_{i}$ through the relationship,
\begin{equation}
	\label{eq:5}
	\Delta Y_{i} = \Delta G_{i} + j \Delta B_{i}, \quad \mbox{where} \quad j = \sqrt{-1}.
\end{equation}

There is also uncertainty associated with the damping $\Delta D_{i}$, the inertia $\Delta H_{i}$ and thus the effective damping $\Delta \beta_{i}$. Finally, there is uncertainty associated with the real and reactive power generation as well as the real and reactive power demand, this introduces uncertainty into the net power driving a generator $\Delta \alpha_{i}$. The vector which represents the uncertainties in the system $\Delta x$ is the following,
\begin{equation}
	\label{eq:7}
	\Delta x = \begin{pmatrix}
		\{\Delta \beta_{i}\} \\
		\{\Delta \alpha_{i}\} \\
		\{\Delta c_{ik}\} \\
		\{\Delta \gamma_{ik}\}
	\end{pmatrix}.
\end{equation}
We denote the uncertainty in the set of power flow state parameters as $\{\Delta y_{i}\}$ so that,
\begin{equation}
	\label{eq:8}
	\{\Delta y_{i}\} = \begin{pmatrix}
		\{\Delta \alpha_{i}\} \\
		\{\Delta c_{ik}\} \\
		\{\Delta \gamma_{ik}\}
	\end{pmatrix}.
\end{equation}
We can then rewrite \cref{eq:7} as
\begin{equation}
	\label{eq:9}
	\Delta x = \begin{pmatrix}
		\{\Delta \beta_{i}\} \\
		\{\Delta y_{i}\}
	\end{pmatrix}.
\end{equation}

We write \cref{eq:9} in this form to explicitly separate the uncertainty associated with the parameters we can control (albeit only up to limited precision), $\{\Delta \beta_{i}\}$ and the parameters that are fixed $\{\Delta y_{i}\}$. As mentioned in \cref{Introduction}, each value in our system is a distribution with associated uncertainty. In the absence of data about uncertainty, we assume that each parameter can be represented by a normal distribution with a mean given by the real power system values and a standard deviation based on the last significant digit of the real power system value. We use the following notation to refer to the normal distribution of the full system $x$,
\begin{equation}
	\label{eq:10}
	x \sim \mathcal{N}(\bar{x},\Delta x),
\end{equation}
with $\bar{x}$ being the mean value and $\Delta x$ being the uncertainty for each element of $x$. In full, $x$ is the following
\begin{equation}
	\label{eq:11}
	x \sim \begin{pmatrix}
		\{\mathcal{N}(\bar{\beta}_{i},\Delta \beta_{i})\} \\
		\{\mathcal{N}(\bar{\alpha}_{i},\Delta \alpha_{i})\} \\
		\{\mathcal{N}(\bar{c}_{ik},\Delta c_{ik})\} \\
		\{\mathcal{N}(\bar{\gamma}_{ik},\Delta \gamma_{ik})\} \\
	\end{pmatrix},
\end{equation}
which is used to test each method's robustness to stability. We generate multiple instances of a targeted steady state with the mean and standard deviation specified and denote each steady state with a superscript $j$. The distributions of the Lyapunov exponent for these instances allow us to determine the robustness of each method for a system under uncertainty.

\section{Optimisation Under Uncertainty}
\label{Optimisation Under Uncertainty}
\subsection{Background}
\label{Background}
In this section, we analyse the impact of uncertainty on optimal power grid configurations. For example, consider the case described in \cref{tab6}. The system parameters that define $\{y_{i}\}$ are given by the mean values measured from the system, but there is actually uncertainty in these values as given in \cref{tab6} for a particular steady state of the power grid. We explore how this uncertainty (and the impossibility of precisely specifying $\beta_{i}$) impacts $\lambda_{L}$ when optimising $\{\beta_{i}\}$. In particular, we look at three levels of uncertainty as specified through the standard deviation, each an order of magnitude apart. Essentially, we are investigating the impact that unwanted deviations in the optimal $\beta_{i}$ as a result of uncertainty have on the overall stability.

\begin{table}[h]
	\footnotesize
	\caption{System model parameters and their corresponding means and uncertainty for a single steady state.}
	\begin{center}
		\begin{tabular}{|c|c|c|}
			\hline
			\textbf{System Parameter} & \textbf{Mean} & \textbf{Uncertainty} \\
			\hline
			Active Power Demand & $160$ & $1$ \\
			\hline
			Reactive Power Demand & $35$ & $1$ \\
			\hline
			Active Power Generation & $630$ & $1$ \\
			\hline
			Reactive Power Generation & $127$ & $1$ \\
			\hline
			Resistance & $0.001$ & $0.001$  \\
			\hline
			Reactance & $0.020$ & $0.001$  \\
			\hline
		\end{tabular}
		\label{tab6}
	\end{center}
\end{table}

We then look at three cases under different levels of uncertainty as summarised in \cref{tab3} which are the $\beta_{=}$, $\beta_{\neq}$ and $\beta_{u}$ cases. In the $\beta_{=}$ case, the same $\beta$ value is used for all generators and this is optimised assuming no uncertainty in $\beta$. In the $\beta_{\neq}$ case, different $\beta$ values are used for each generator and this is also optimised assuming no uncertainty. In the $\beta_{u}$ case, different distributions of $\beta$ are used for each generator and the mean of each $\beta$ distribution $\langle \beta\rangle$ is optimised assuming some uncertainty in $\beta$. All cases are tested by calculating $\lambda_{L}$ for a specified noise in $\{\beta_{i}\}$ and $\{y_{i}\}$ for a given steady state $j$ and running the test $10000$ times to generate distributions of $\lambda_{L}$. All calculations were performed in Matlab and the code is provided at \url{https://github.com/John-Moloney/Power-Grid-Stability-Under-Uncertainty}.

It is a challenge to compare these distributions as there are competing factors that determine whether a distribution corresponds to an optimal physical state, two plausible factors are the distribution that produces the most optimal $\lambda_{L}$ value and the distribution that is least likely to produce the worst $\lambda_{L}$. We introduce quantile plots as a way of comparing performance when distinct methods are used or different numbers of samples are taken. Quantile plots give us the $\lambda_{L}$ value for a given quantile from a $\lambda_{L}$ distribution. The significance of these plots is that they show us the methods that have the best performance as given by the curves with the lowest $\lambda_{L}$ and they also help us identify under what circumstances these methods perform well. For example, if a method performs well at high quantiles then this method performs better than other methods even at its worst values and we would describe this method as robust. However, if a method performs well at low quantiles then this method performs better close to its theoretical optimum than other methods near their theoretical optima.

\begin{table}[h]
	\footnotesize
	\caption{The three cases that are tested for their stability properties are; a single $\beta_{i}$ value applied to each generator in the system denoted as the $\beta_{=}$ case, an optimal set of $\beta_{i}$ values generated from a system with no uncertainty denoted as the $\beta_{\neq}$ case and an optimal set of $\{\langle \beta_{i}\rangle\}$ generated from a system with uncertainty denoted as $\beta_{u}$.}
	\begin{center}
		\begin{tabular}{|c|c|c|c|}
			\hline
			\textbf{Case} & \textbf{Know} & \textbf{Optimise} & \textbf{Test} \\
			\hline
			$\beta_{=}$ & $x^{j} \sim \mathcal{N}(\bar{x}^{j},0)$ & $\beta_{i}$ & $\lambda_{L}(\{x^{j} + \Delta x^{j}\})$ \\
			\hline
			$\beta_{\neq}$ & $x^{j} \sim \mathcal{N}(\bar{x}^{j},0)$ & $\{\beta_{i}\}$ & $\lambda_{L}(\{x^{j} + \Delta x^{j}\})$ \\
			\hline
			 $\beta_{u}$ & $x^{j} \sim \mathcal{N}(\bar{x}^{j},\Delta x^{j})$ & $\{\langle \beta_{i}\rangle\}$ \mbox{given} $\{\Delta \beta_{i}\}$  & $\lambda_{L}(\{x^{j} + \Delta x^{j}\})$ \\
			\hline
		\end{tabular}
		\label{tab3}
	\end{center}
\end{table}

\subsection{Base Cases}
\label{Base Cases}
We begin by considering methods that assume zero uncertainty when calculating optimal $\beta$ values. We look at how these methods perform in the presence of uncertainty. We do this by calculating the distributions of $\lambda_{L}$ for both the $\mathcal{N}(\bar{\beta}_{=},0)$ and $\mathcal{N}(\bar{\beta}_{\neq},0)$ cases by specifying the uncertainty in each $x^{j}$, $\mathcal{N}(\bar{x}^{j},\Delta x^{j})$ and calculating $\lambda_{L}$ $10000$ times.

\begin{figure}[h]
	\centering
	\begin{tikzpicture}[scale=1.0]
		
		\centering
		
		\begin{axis}[ymax=3600,ymin=0, xmax=-1.2, xmin=-4.4, xlabel={$\lambda_{L}$}, ylabel={Frequency}, ymajorgrids=false, xmajorgrids=false, xtick pos=bottom, ytick pos=left, xtick = {-4.4, -4, -3.6, -3.2, -2.8, -2.4, -2, -1.6},scatter/classes={a={mark=square,red},b={mark=o,blue},c={mark=x,black}}]
			\addplot+[only marks,scatter,scatter src=explicit symbolic,mark=square] coordinates { (-3.885, 1207)[a] (-3.855, 3575)[a] (-3.825, 58)[a] (-3.795, 120)[a] (-3.765, 129)[a] (-3.735, 164)[a] (-3.705, 218)[a] (-3.675,248)[a] (-3.645,308)[a] (-3.615,300)[a] (-3.585,332)[a] (-3.555,375)[a] (-3.525,354)[a] (-3.495,343)[a] (-3.465, 377)[a] (-3.435, 339)[a] (-3.405, 345)[a] (-3.375,254)[a] (-3.345,251)[a] (-3.315,237)[a] (-3.285,135)[a] (-3.255,126)[a] (-3.225,79)[a] (-3.195,64)[a] (-3.165,39)[a] (-3.135,8)[a] (-3.105,11)[a] (-3.3075,3)[a] (-3.3045,1)[a] (-3.015,14)[a]};
			
			\addplot+[only marks,scatter,scatter src=explicit symbolic,mark=circle] coordinates { (-3.885, 67)[b] (-3.855, 1232)[b] (-3.825, 2534)[b] (-3.795, 1039)[b] (-3.765, 241)[b] (-3.735, 155)[b] (-3.705, 153)[b] (-3.675,178)[b] (-3.645,200)[b] (-3.615,230)[b] (-3.585,258)[b] (-3.555,286)[b] (-3.525,310)[b] (-3.495,323)[b] (-3.465, 290)[b] (-3.435, 325)[b] (-3.405, 315)[b] (-3.375,307)[b] (-3.345,268)[b] (-3.315,265)[b] (-3.285,229)[b] (-3.255,210)[b] (-3.225,173)[b] (-3.195,115)[b] (-3.165,108)[b] (-3.135,57)[b] (-3.105,55)[b] (-3.3075,33)[b] (-3.3045,17)[b] (-3.015,14)[b] (-2.985,6)[b] (-2.955,3)[b]};
			
			\addplot+[only marks,scatter,scatter src=explicit symbolic,mark=triangle] coordinates { (-3.875,3)[c] (-3.825,29)[c] (-3.775,51)[c] (-3.725,98)[c] (-3.675,148)[c] (-3.625,222)[c] (-3.575,267)[c] (-3.525,364)[c] (-3.475,430)[c] (-3.425,475)[c] (-3.375,527)[c] (-3.325,585)[c] (-3.275,540)[c] (-3.225,563)[c] (-3.175,479)[c] (-3.125,489)[c] (-3.075,469)[c] (-3.025,434)[c] (-2.975,420)[c] (-2.925,380)[c] (-2.875,359)[c] (-2.825,350)[c] (-2.775,334)[c] (-2.725,315)[c] (-2.675,292)[c] (-2.625,255)[c] (-2.575,230)[c] (-2.525,218)[c] (-2.475,177)[c] (-2.425,159)[c] (-2.375,127)[c] (-2.325,80)[c] (-2.275,54)[c] (-2.225,35)[c] (-2.175,22)[c] (-2.125,11)[c] (-2.075,6)[c] (-2.025,2)[c] (-1.975,1)[c]};
			\legend{$\sigma = 0.01$,$\sigma = 0.1$,$\sigma = 1$}
		\end{axis}	
		
	\end{tikzpicture}
	
	\caption{Comparing the distribution of $\lambda_{L}$ values for $\beta_{=}$ with $\sigma = 0.01$, $\sigma = 0.1$ and $\sigma = 1$.}
	\label{fig:11}
\end{figure}

\begin{figure}[h]
	\centering
	\begin{tikzpicture}[scale=1.0]
		
		\centering
		
		\begin{axis}[ymax=800,ymin=0, xmax=-1.2, xmin=-4.4, xlabel={$\lambda_{L}$}, ylabel={Frequency}, ymajorgrids=false, xmajorgrids=false, xtick pos=bottom, ytick pos=left, xtick={-4.4,-4,-3.6,-3.2,-2.8,-2.4,-2,-1.6,-1.2},scatter/classes={a={mark=square,red},b={mark=o,blue},c={mark=x,black}}]
			\addplot+[only marks,scatter,scatter src=explicit symbolic] coordinates { (-4.275,15)[a] (-4.225,318)[a] (-4.175,522)[a] (-4.125,417)[a] (-4.075,486)[a] (-4.025,574)[a] (-3.975,691)[a] (-3.925,704)[a] (-3.875,743)[a] (-3.825,685)[a] (-3.775,609)[a] (-3.725,530)[a] (-3.675,480)[a] (-3.625,501)[a] (-3.575,496)[a] (-3.525,557)[a] (-3.475,488)[a] (-3.425,363)[a] (-3.375,317)[a] (-3.325,219)[a] (-3.275,141)[a] (-3.225,81)[a] (-3.175,45)[a] (-3.125,11)[a] (-3.075,3)[a] (-3.025,3)[a] (-2.975,1)[a]};
			
			\addplot+[only marks,scatter,scatter src=explicit symbolic] coordinates { (-4.275,15)[b] (-4.225,82)[b] (-4.175,238)[b] (-4.125,351)[b] (-4.075,380)[b] (-4.025,426)[b] (-3.975,549)[b] (-3.925,641)[b] (-3.875,698)[b] (-3.825,724)[b] (-3.775,680)[b] (-3.725,601)[b] (-3.675,516)[b] (-3.625,475)[b] (-3.575,419)[b] (-3.525,401)[b] (-3.475,429)[b] (-3.425,389)[b] (-3.375,383)[b] (-3.325,394)[b] (-3.275,331)[b] (-3.225,295)[b] (-3.175,227)[b] (-3.125,164)[b] (-3.075,97)[b] (-3.025,50)[b] (-2.975,25)[b] (-2.925,14)[b] (-2.875,4)[b] (-2.825,2)[b]};
			
			\addplot+[only marks,scatter,scatter src=explicit symbolic] coordinates { (-4.25,1)[c] (-4.15,12)[c] (-4.05,9)[c] (-3.95,44)[c] (-3.85,105)[c] (-3.75,236)[c] (-3.65,408)[c] (-3.55,544)[c] (-3.45,637)[c] (-3.35,684)[c] (-3.25,697)[c] (-3.15,729)[c] (-3.05,714)[c] (-2.95,668)[c] (-2.85,674)[c] (-2.75,663)[c] (-2.65,661)[c] (-2.55,554)[c] (-2.45,527)[c] (-2.35,459)[c] (-2.25,396)[c] (-2.15,312)[c] (-2.05,159)[c] (-1.95,75)[c] (-1.85,24)[c] (-1.75,4)[c] (-1.65,1)[c] (-1.55,2)[c] (-1.35,1)[c]};
			\legend{$\sigma = 0.01$,$\sigma = 0.1$,$\sigma = 1$}
		\end{axis}	
		
	\end{tikzpicture}
	
	\caption{Comparing the distribution of $\lambda_{L}$ values for $\beta_{\neq}$ with $\sigma = 0.01$, $\sigma = 0.1$ and $\sigma = 1$.}
	\label{fig:12}
\end{figure}

\Cref{fig:11} and \cref{fig:12} tests both $\mathcal{N}(\bar{\beta}_{=},0)$ and $\mathcal{N}(\bar{\beta}_{\neq},0)$ by applying noise to each $\beta_{i}$ value with a mean of $\bar{\beta}_{=}$ and $\bar{\beta}_{\neq}$ respectively with standard deviations of $\sigma = 0.01$, $\sigma = 0.1$ and $\sigma = 1$ as well as noise in each $y_{i}^{j}$ according to \cref{tab6}. It is clear from \cref{fig:11} that for small amounts of noise applied to $\beta_{=}$ ($\sigma = 0.01$), there is a concentration of $\lambda_{L}$ around the minimum $\lambda_{L}$ value. However, as the applied noise in the damping coefficients increases, we see for $\sigma = 0.1$, the $\lambda_{L}$ distribution is qualitatively similar but shifted to the right. The case where the applied noise is $\sigma = 1$ has a much more Gaussian-like distribution with the mean shifted to the right of the other distributions with smaller deviations. \Cref{fig:12} follows a similar trend to that of \cref{fig:11} as the distributions shift to the right as the uncertainty increases. However, instead of a concentration of $\lambda_{L}$ resulting in a spike at the left-hand side of the distribution around the most optimal $\lambda_{L}$ values, there are instead more normal-like distributions for $\sigma = 0.01$ and $\sigma = 0.1$.

\subsection{Simulated Annealing Under Uncertainty Algorithm}
\label{Uncertainty in power flow}
This motivates the idea that there may be a better choice of $\{\beta_{i}\}$ that accounts for the fact that there is uncertainty in the system. Therefore we look at using an optimisation algorithm that optimises the choice of $\{\beta_{i}\}$ under uncertainty. We use a simulated annealing algorithm developed by Tchechmedjiev et al. \cite{Tchechmedjiev2012}, the pseudocode of which is provided in \cref{alg:optimisation}. The main differences between this algorithm and a standard simulated annealing algorithm are that an extra loop is introduced to sample the objective function.

\begin{algorithm}
	\caption{Simulated Annealing Under Uncertainty}
	\label{alg:optimisation}
	\begin{algorithmic}[1]
		\STATE{Let $\{\beta_i\} = \{\beta_{0}\}$} 
		\FOR{$k = 0$ through $k_{max}$}
		\STATE{$T\gets$ temperature$(\frac{b_{0} 2 \sigma_{0}}{k^{T} \sqrt{N}})$}
		\STATE{Pick a random neighbour, $\{\beta_{i}\}\gets$ $\{\beta_{\text{new}}\}$}
		\FOR{$n = 1$ through $n_{max}$}
		\STATE{$\lambda_{L}\gets$ $\lambda(\{\beta_{\text{new}}\})$}
		\ENDFOR\label{Sampling}
		
		$\lambda_{\text{new}}\gets$ mean($\lambda_{L}$)
		
		If $P(E(\{\beta_{i}\}),\{\beta_{\text{new}}\},T)\geq$ random($0,1$)
		
		$\lambda\gets\lambda_{\text{new}}$
		\ENDFOR\label{Cooling Schedule}
		\STATE{\textbf{return} $\{\beta_{i}\}$}
	\end{algorithmic}
\end{algorithm}

\begin{figure}[h!]
	\centering
	\vspace{0mm}
	\includegraphics[width=0.75\textwidth]{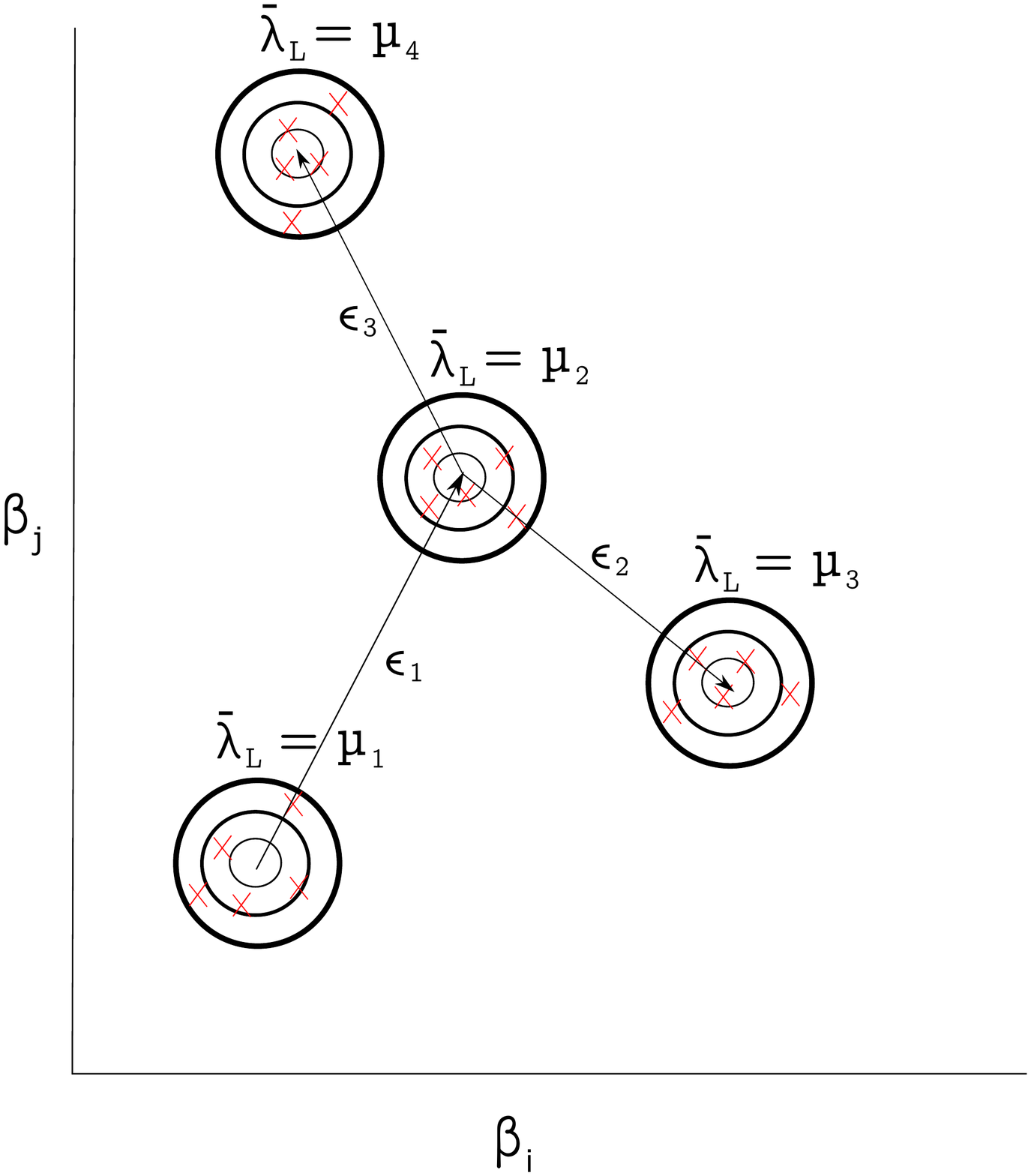}
	\caption{Example of a two generator system with coordinates given in the $(\beta_{i},\beta_{j})$ space. The Gaussian distribution centred at $(\beta_{i},\beta_{j})$ is sampled many times with each sample represented by a red $x$, the average $\lambda_{L}$ is then taken as given by $\bar{\lambda}_{L}$. We then perturb this point by $\epsilon_{i}$ and we once again sample a Gaussian distribution at this new point and calculate the new $\bar{\lambda}_{L}$. If the new $\bar{\lambda}_{L}$ is less than the old $\bar{\lambda}_{L}$ then we accept the new $\bar{\lambda}_{L}$ as is done in the diagram when perturbing by $\epsilon_{1}$. However, if $\bar{\lambda}_{L}$ is more than the old $\bar{\lambda}_{L}$ then we may still accept the new $\bar{\lambda}_{L}$ with some probability which depends on the state of the simulated annealing algorithm. This happens in the diagram when comparing $\mu_{4}$ and $\mu_{2}$ unlike the case that compares $\mu_{3}$ and $\mu_{2}$ assuming $\mu_{4} > \mu_{3} > \mu_{1} > \mu_{2}$.}
	\label{fig:Mean_Gaussian_SA}
\end{figure}

This is best illustrated by \cref{fig:Mean_Gaussian_SA} which gives an example of how the algorithm works for the first three iterations. Some initial set of $\{\beta_{i}\}$ is chosen as the set of mean values of the Gaussian distributions for $\{\beta_{i}\}$. The distribution of each $\beta_{i}$ is sampled and for the set of sampled $\{\beta_{i}\}$, a $\lambda_{L}$ value is calculated. This is repeated many times and the average $\lambda_{L}$ is recorded which is denoted as $\bar{\lambda}_{L}$. There is an $\epsilon$ value imposed for each $\beta_{i}$ that determines the perturbation to the corresponding effective damping parameter. We repeat the sampling process for the new set of $\beta_{i}$ following the perturbation and compare the new and old $\bar{\lambda}_{L}$ values to decide whether to accept or reject the new set of $\beta_{i}$. We decide whether to accept the new $\bar{\lambda}_{L}$ using the metropolis criterion and Boltzmann distribution,
\begin{equation}
	\label{eq:18}
	P = \begin{cases}
		1 & \mbox{if} \quad \Delta \lambda < 0\\
		\exp{\frac{-\Delta \lambda}{T}} & \mbox{otherwise},
	\end{cases}
\end{equation}
where $\Delta \lambda$ is the difference between the new and old $\bar{\lambda}_{L}$. We then repeat this process until the stopping conditions of the simulated annealing algorithm are satisfied.

As we are dealing with an uncertain objective function, we must sample the objective function to determine the expected value which we use as a global objective function. We use the global objective term $\Phi_{E}$ in the simulated annealing algorithm to determine whether to accept the new $\bar{\lambda}_{L}$. This $\Phi_{E}$ is defined with an expected value and a penalty term, the penalty term is added to deal with the variability of the distributions. The expected value is given by,
\begin{equation}
	\label{eq:14}
	\bar{r}_{0} = \frac{1}{N} \sum_{r_{i}\in L}^{}r_{i}.
\end{equation}
with the set $L$ consisting of $N$ samples generated when evaluating $\lambda_{L}$,
\begin{equation}
	\label{eq:15}
	L = \cup_{i = 1}^{N}\{\lambda_{L}\}.
\end{equation}
The penalty term $\Psi_{E}$ is given by,
\begin{equation}
	\label{eq:12}
	\Psi_{E} = \frac{b_{0}}{k^{T}} \frac{2 \sigma_{0}}{\sqrt{N}},
\end{equation}
where $b_{0}$ is a small constant, k is a constant that represents the rate of increase of $b(T)$, $T$ is the temperature of the system and the scaled standard deviation is,
\begin{equation}
	\label{eq:13}
	\sigma_{0} = \frac{\sqrt{\sum_{r_{i}\in L}^{} (r_{i} - \bar{r}_{0})^{2}}}{N}.
\end{equation}
The expected value and penalty term combine to give the global objective function,
\begin{equation}
	\label{eq:17}
	\Phi_{E} = \bar{r}_{0} + \Psi_{E}.
\end{equation}

There are several ways to pick the perturbation of each $\epsilon$, one such way is to pick a random point on a hypersphere. The benefit of picking a random point on a hypersphere is that this is an efficient way of picking a random perturbation for each $\beta_{i}$. The step may be positive or negative and the normalisation ensures the steps are not too small or too large while also ensuring there is no significant overlap between subsequent picking of distributions. A simple method for selecting points on a $4$ dimensional sphere is described by Marsaglia \cite{Marsaglia1972} but to pick a random point on a $n$ dimensional hypersphere, it is necessary to extend this method. We first generate $n$ random Gaussian variables $a_{1}, a_{2}, ... a_{n}$ and then the distributions of the following vectors,
\begin{equation}
	\label{eq:16}
	\epsilon = \frac{1}{\sqrt{a_{1}^{2}, a_{2}^{2}, ... a_{n}^{2}}}\begin{bmatrix}
		a_{1} \\
		a_{2} \\
		\vdots \\
		a_{n}
	\end{bmatrix}
\end{equation}
are the points we use to perturb $x^{j}$.

One drawback to this method is that it can quickly become computationally expensive if we wish to get as good an approximation of the global minimum as possible. This is due to the sparsity of data whereby as the number of dimensions increases, the volume of the model space which must be explored increases while the number of optimal solutions does not increase \cite{Curtis2001}. This leads to the desired data becoming more sparse which makes it a difficult task to approximate the global optimum. More samples are then required to accurately estimate the summary statistic which is the mean in this instance.

\begin{figure}[h]
	\centering
	\begin{tikzpicture}
		\centering
		\begin{axis}[
			legend style={nodes={scale=0.5, transform shape}},
			legend image post style={scale=0.5},
			scaled ticks=false, 
			tick label style={/pgf/number format/fixed},
			xlabel={Quantile},
			ylabel={$\lambda_{L}$},
			xmin = 10, xmax = 90,
			ymin = -3.6, ymax = -2.2,
			legend pos=north west,
			]	
			\addplot[line width=1pt,solid,color=blue,mark=square*,dashed,mark options=solid]%
			coordinates {(10, -3.5688) (20, -3.4597) (30, -3.3687) (40, -3.2811) (50, -3.1952) (60, -3.1006) (70, -3.001) (80, -2.8734) (90, -2.6848)};
			\addplot[line width=1pt,solid,color=red,mark=square*,solid,mark options=solid]%
			coordinates {(10, -3.4822) (20, -3.365) (30, -3.2589) (40, -3.1622) (50, -3.0596) (60, -2.9534) (70, -2.8295) (80, -2.683) (90, -2.4901)};
			\addplot[line width=1pt,solid,color=blue,mark=square*,dashed,mark options=solid]%
			coordinates {(10, -3.4797) (20, -3.3511) (30, -3.242) (40, -3.134) (50, -3.0306) (60, -2.9112) (70, -2.7864) (80, -2.6295) (90, -2.426)};
			\addplot[line width=1pt,solid,color=blue,mark=square*,dashed,mark options=solid]%
			coordinates {(10, -3.572) (20, -3.4635) (30, -3.3559) (40, -3.2501) (50, -3.1391) (60, -3.0099) (70, -2.8601) (80, -2.6973) (90, -2.4868)};
			\addplot[line width=1pt,solid,color=blue,mark=square*,dashed,mark options=solid]%
			coordinates {(10, -3.5045) (20, -3.3969) (30, -3.2944) (40, -3.1938) (50, -3.0979) (60, -2.9893) (70, -2.8638) (80, -2.7223) (90, -2.532)};
			\addplot[line width=1pt,solid,color=blue,mark=square*,dashed,mark options=solid]%
			coordinates {(10, -3.5125) (20, -3.3831) (30, -3.2787) (40, -3.1766) (50, -3.0755) (60, -2.9592) (70, -2.8348) (80, -2.6892) (90, -2.5011)};
			\addplot[line width=1pt,solid,color=red,mark=square*,solid,mark options=solid]%
			coordinates {(10, -3.5404) (20, -3.4202) (30, -3.3161) (40, -3.2096) (50, -3.1002) (60, -2.9862) (70, -2.8613) (80, -2.7045) (90, -2.5088)};
			\addplot[line width=1pt,solid,color=red,mark=square*,solid,mark options=solid]%
			coordinates {(10, -3.4495) (20, -3.3321) (30, -3.2254) (40, -3.1244) (50, -3.0235) (60, -2.9171) (70, -2.8012) (80, -2.6507) (90, -2.4677)};
			\addplot[line width=1pt,solid,color=red,mark=square*,solid,mark options=solid]%
			coordinates {(10, -3.5459) (20, -3.4242) (30, -3.3186) (40, -3.2187) (50, -3.1232) (60, -3.0089) (70, -2.88) (80, -2.7328) (90, -2.5253)};
			\addplot[line width=1pt,solid,color=red,mark=square*,solid,mark options=solid]%
			coordinates {(10, -3.228) (20, -3.0873) (30, -2.9691) (40, -2.867) (50, -2.77) (60, -2.6714) (70, -2.5716) (80, -2.458) (90, -2.2917)};
			\addlegendentry{$n = 100$}
			\addlegendentry{$n = 20$}
		\end{axis}
	\end{tikzpicture}
	\caption{$\lambda_{L}$ versus quantile plot for perturbations about a given $x^{j}$ with $\sigma = 1$ for the $\beta_{u}$ case with each $\epsilon$ picked from a Gaussian distribution with the $\bar{\lambda}_{L}$ calculated from sampling $20$ times for the solid red curves and $100$ times for the dashed blue curves.}
	\label{fig:10}
\end{figure}
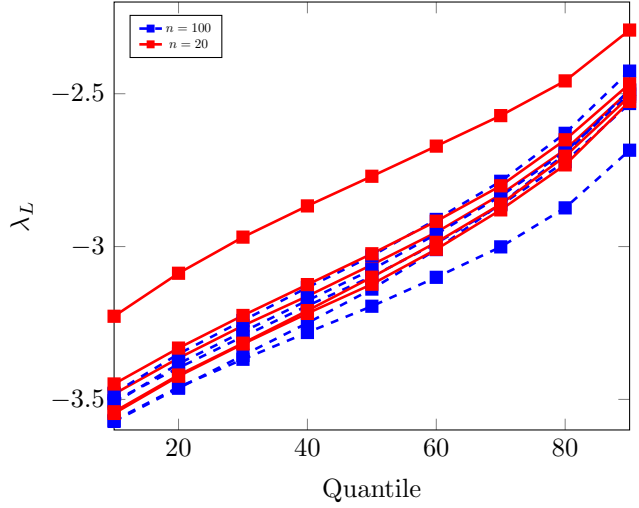

The benefit of taking more samples is clear from the quantile plot in \cref{fig:10} which compares sampling $\lambda_{L}$ $20$ times before calculating $\bar{\lambda}_{L}$ as given by the red curves to sampling $\lambda_{L}$ $100$ times before calculating $\bar{\lambda}_{L}$ as given by the blue curves. The $20$ samples case produces the worst result for a run and the $100$ samples case produces the best result of any run however, there is little difference between sampling $20$ and $100$ times for the majority of cases. This shows that while it is possible to achieve better results by using more samples, this is difficult given the sparsity of data for large dimensions. It would therefore be better to combine large numbers of samples with many runs to achieve the best results. This would not be feasible as the computational resources needed would scale with both the number of samples and the number of runs which would be an inefficient and computationally expensive process. In this paper, we use 100 samples when calculating the optimal $\beta_{i}$ under uncertainty although more samples could be taken if a more robust stability condition is desired and the necessary computational resources are available.

\subsection{Uncertainty in Effective Damping}
\label{Uncertainty in Effective Damping}
We use this simulated annealing under uncertainty algorithm to determine the best centre of each $\beta_{i}$ distribution to get the best $\lambda_{L}$ under uncertainty. In this case, the power flow state parameters are known with a high degree of certainty and it is $\{\beta_{i}\}$ that has noise which impacts stability, so that $y_{i}$ is drawn from $\mathcal{N}(\{\bar{y}_{i}\},0)$ and $\beta_{i}$ is drawn from $\mathcal{N}(\{\bar{\beta}_{i}\}, \Delta \beta_{i})$.

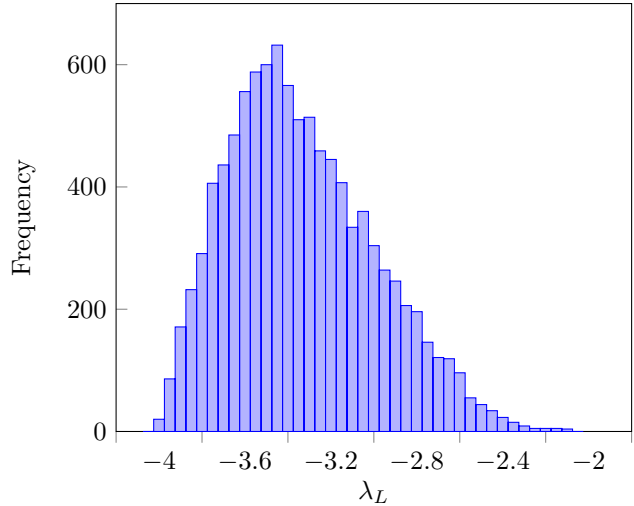
\begin{figure}[h]
	\centering
	\begin{tikzpicture}[scale=1.0]
		
		\centering
		
		\begin{axis}[ybar interval, ymax=700,ymin=0, xmax=-1.6, xmin=-4, xlabel={$\lambda_{L}$}, ylabel={Frequency}, ymajorgrids=false, xmajorgrids=false, xtick pos=bottom, ytick pos=left, xtick = {-4, -3.6, -3.2, -2.8, -2.4, -2, -1.6}]
			\addplot coordinates { (-3.875, 0) (-3.825, 20) (-3.775, 86) (-3.725, 171) (-3.675, 232) (-3.625, 291) (-3.575, 406) (-3.525,436) (-3.475,485) (-3.425,556) (-3.375,588) (-3.325,600) (-3.275,632) (-3.225,566) (-3.175, 510) (-3.125, 514) (-3.075, 459) (-3.025,445) (-2.975,407) (-2.925,334) (-2.875,360) (-2.825,304) (-2.775,264) (-2.725,246) (-2.675,206) (-2.625,196) (-2.575,146) (-2.525,121) (-2.475,119) (-2.425,96) (-2.375,55) (-2.325,44) (-2.275,34) (-2.225,23) (-2.175,15) (-2.125,9) (-2.075,5) (-2.025,5) (-1.975,5) (-1.925,4) (-1.875,0) (-1.825,1)};
		\end{axis}
		
	\end{tikzpicture}
	
	\caption{Histogram of $\lambda_{L}$ values when perturbing $x$ with $\sigma = 1$ for the distribution of $\{\beta_{i}\}$ after optimising $\{\beta_{i}\}$ for $\beta_{u}$ case.}
	\label{fig:6}
\end{figure}

It is clear from \cref{fig:6} that this is also a positively skewed distribution with $\lambda_{L}$ values more likely to occur closer to the lowest $\lambda_{L}$ value as opposed to the highest. It is also clear that the mode occurs at approximately $\lambda_{L} = -3.5$. The centres of the $\{\beta_{i}\}$ values for the optimal configuration associated with each method are shown in \cref{tab1}, \cref{tab10} and \cref{tab4} for $\sigma = 1$, $\sigma = 0.1$ and $\sigma = 0.01$ respectively. We analyse the quantitative differences between $\lambda_{L}$ distributions for the $\beta_{=}$, $\beta_{\neq}$ and $\beta_{u}$ methods when we conduct statistical analysis in \cref{Statistical Analysis}.

\begin{table}[h]
	\footnotesize
	\caption{The $\beta$ value for each generator in the $\beta_{=}$, $\beta_{\neq}$ and $\beta_{u}$ case for $\sigma = 1$ when $\epsilon$ is picked from a normal distribution.}
	\begin{center}
		\begin{tabular}{|c|c|c|c|}
			\hline
			\textbf{Generator} & \textbf{\textit{$\beta_{=}$}} & \textbf{\textit{$\beta_{\neq}$}} & \textbf{\textit{$\beta_{u}$}} \\
			\hline
			$1$ & $7.75$ & $12.83$ & $7.02$ \\
			\hline
			$2$ & $7.75$ & $13.54$ & $10.60$ \\
			\hline
			$3$ & $7.75$ & $9.76$ & $9.37$ \\
			\hline
			$4$ & $7.75$ & $12.64$ & $13.37$ \\
			\hline
			$5$ & $7.75$ & $12.23$ & $7.83$ \\
			\hline
			$6$ & $7.75$ & $7.01$  & $6.73$ \\
			\hline
			$7$ & $7.75$ & $10.38$  & $8.99$ \\
			\hline
			$8$ & $7.75$ & $8.03$ & $10.32$ \\
			\hline
			$9$ & $7.75$ & $8.29$ & $8.84$ \\
			\hline
			$10$ & $7.75$ & $5.89$ & $6.06$ \\
			\hline
		\end{tabular}
		\label{tab1}
	\end{center}
\end{table}

\begin{table}[h]
	\footnotesize
	\caption{The $\beta$ value for each generator in the $\beta_{=}$, $\beta_{\neq}$ and $\beta_{u}$ case for $\sigma = 0.1$ when $\epsilon$ is picked from a normal distribution.}
	\begin{center}
		\begin{tabular}{|c|c|c|c|}
			\hline
			\textbf{Generator} & \textbf{\textit{$\beta_{=}$}} & \textbf{\textit{$\beta_{\neq}$}} & \textbf{\textit{$\beta_{u}$}} \\
			\hline
			$1$ & $7.75$ & $12.83$ & $8.04$ \\
			\hline
			$2$ & $7.75$ & $13.54$ & $10.52$ \\
			\hline
			$3$ & $7.75$ & $9.76$ & $9.22$ \\
			\hline
			$4$ & $7.75$ & $12.64$ & $8.36$ \\
			\hline
			$5$ & $7.75$ & $12.23$ & $8.02$ \\
			\hline
			$6$ & $7.75$ & $7.01$  & $8.11$ \\
			\hline
			$7$ & $7.75$ & $10.38$  & $7.70$ \\
			\hline
			$8$ & $7.75$ & $8.03$ & $10.16$ \\
			\hline
			$9$ & $7.75$ & $8.29$ & $7.38$ \\
			\hline
			$10$ & $7.75$ & $5.89$ & $7.20$\\
			\hline
		\end{tabular}
		\label{tab10}
	\end{center}
\end{table}

\begin{table}[h]
	\footnotesize
	\caption{The $\beta$ value for each generator in the $\beta_{=}$, $\beta_{\neq}$ and $\beta_{u}$ case for $\sigma = 0.01$ when $\epsilon$ is picked from a normal distribution.}
	\begin{center}
		\begin{tabular}{|c|c|c|c|}
			\hline
			\textbf{Generator} & \textbf{\textit{$\beta_{=}$}} & \textbf{\textit{$\beta_{\neq}$}} & \textbf{\textit{$\beta_{u}$}} \\
			\hline
			$1$ & $7.75$ & $12.83$ & $8.02$ \\
			\hline
			$2$ & $7.75$ & $13.54$ & $14.34$ \\
			\hline
			$3$ & $7.75$ & $9.76$ & $12.02$ \\
			\hline
			$4$ & $7.75$ & $12.64$ & $8.23$ \\
			\hline
			$5$ & $7.75$ & $12.23$ & $7.24$ \\
			\hline
			$6$ & $7.75$ & $7.01$  & $8.86$ \\
			\hline
			$7$ & $7.75$ & $10.38$  & $9.19$ \\
			\hline
			$8$ & $7.75$ & $8.03$ & $8.63$ \\
			\hline
			$9$ & $7.75$ & $8.29$ & $8.07$ \\
			\hline
			$10$ & $7.75$ & $5.89$ & $6.81$ \\
			\hline
		\end{tabular}
		\label{tab4}
	\end{center}
\end{table}

\section{Statistical Analysis}
\label{Statistical Analysis}
\begin{figure}[h]
	\centering
	\begin{tikzpicture}[scale=0.5]
		\centering
		\begin{axis}[
			legend style={nodes={scale=0.5, transform shape}}, 
			legend image post style={scale=0.5},
			scaled ticks=false, 
			tick label style={/pgf/number format/fixed},
			xlabel={Quantile},
			ylabel={$\lambda_{L}$},
			xmin = 10, xmax = 90,
			ymin = -3.6, ymax = -2.2,
			legend pos=north west,
			]	
			\node at (5mm,52.5mm) {\Large\bfseries\sffamily\selectfont a};
			\addplot[line width=1pt,solid,color=brown,mark=triangle*,solid,mark options=solid]%
			coordinates {(10, -3.57) (20, -3.4) (30, -3.25) (40, -3.11) (50, -2.98) (60, -2.83) (70, -2.67) (80, -2.5) (90, -2.31)};
			\addplot[line width=1pt,solid,color=blue,mark=square*,dashed,mark options=solid]%
			coordinates {(10, -3.5235) (20, -3.4068) (30, -3.3127) (40, -3.2254) (50, -3.1278) (60, -3.0247) (70, -2.9105) (80, -2.7628) (90, -2.5872)};
			\addplot[line width=1pt,solid,color=red,mark=*,dotted,mark options=solid]%
			coordinates {(10, -3.5436) (20, -3.4243) (30, -3.3299) (40, -3.2361) (50, -3.1378) (60, -3.0349) (70, -2.9157) (80, -2.7743) (90, -2.5787)};
		\end{axis}
	\end{tikzpicture}
	\begin{tikzpicture}[scale=0.5]
		\centering
		\begin{axis}[
			legend style={nodes={scale=0.5, transform shape}}, 
			legend image post style={scale=0.5},
			scaled ticks=false, 
			tick label style={/pgf/number format/fixed},
			xlabel={Quantile},
			ylabel={$\lambda_{L}$},
			xmin = 10, xmax = 90,
			ymin = -3.6, ymax = -2.2,
			legend pos=north west,
			]	
			\node at (5mm,52.5mm) {\Large\bfseries\sffamily\selectfont b};
			\addplot[line width=1pt,solid,color=brown,mark=triangle*,solid,mark options=solid]%
			coordinates {(10, -3.57) (20, -3.4) (30, -3.25) (40, -3.11) (50, -2.98) (60, -2.83) (70, -2.67) (80, -2.5) (90, -2.31)};
			\addplot[line width=1pt,solid,color=blue,mark=square*,dashed,mark options=solid]%
			coordinates {(10, -3.5195) (20, -3.4083) (30, -3.3152) (40, -3.223) (50, -3.1244) (60, -3.0169) (70, -2.894) (80, -2.7535) (90, -2.5848)};
			\addplot[line width=1pt,solid,color=red,mark=*,dotted,mark options=solid]%
			coordinates {(10, -3.6163) (20, -3.5004) (30, -3.3998) (40, -3.2915) (50, -3.1744) (60, -3.0328) (70, -2.8751) (80, -2.7121) (90, -2.534)};
		\end{axis}
	\end{tikzpicture}
	\begin{tikzpicture}[scale=0.5]
		\centering
		\begin{axis}[
			legend style={nodes={scale=0.5, transform shape}}, 
			legend image post style={scale=0.5},
			scaled ticks=false, 
			tick label style={/pgf/number format/fixed},
			xlabel={Quantile},
			ylabel={$\lambda_{L}$},
			xmin = 10, xmax = 90,
			ymin = -3.6, ymax = -2.2,
			legend pos=north west,
			]	
			\node at (5mm,52.5mm) {\Large\bfseries\sffamily\selectfont c};
			\addplot[line width=1pt,solid,color=brown,mark=triangle*,solid,mark options=solid]%
			coordinates {(10, -3.57) (20, -3.4) (30, -3.25) (40, -3.11) (50, -2.98) (60, -2.83) (70, -2.67) (80, -2.5) (90, -2.31)};
			\addplot[line width=1pt,solid,color=blue,mark=square*,dashed,mark options=solid]%
			coordinates {(10, -3.5195) (20, -3.4083) (30, -3.3152) (40, -3.223) (50, -3.1244) (60, -3.0169) (70, -2.894) (80, -2.7535) (90, -2.5848)};
			\addplot[line width=1pt,solid,color=red,mark=*,dotted,mark options=solid]%
			coordinates {(10, -3.5826) (20, -3.4564) (30, -3.3371) (40, -3.2148) (50, -3.0782) (60, -2.9181) (70, -2.7525) (80, -2.5857) (90, -2.4005)};
		\end{axis}
	\end{tikzpicture}
	\begin{tikzpicture}[scale=0.5]
		\centering
		\begin{axis}[
			legend style={nodes={scale=0.5, transform shape}}, 
			legend image post style={scale=0.5},
			scaled ticks=false, 
			tick label style={/pgf/number format/fixed},
			xlabel={Quantile},
			ylabel={$\lambda_{L}$},
			xmin = 10, xmax = 90,
			ymin = -4.2, ymax = -3.2,
			legend pos=north west,
			]	
			\node at (5mm,52.5mm) {\Large\bfseries\sffamily\selectfont d};
			\addplot[line width=1pt,solid,color=brown,mark=triangle*,solid,mark options=solid]%
			coordinates {(10, -4.0595) (20, -3.9662) (30, -3.8897) (40, -3.8194) (50, -3.7413) (60, -3.6472) (70, -3.5238) (80, -3.4019) (90, -3.2751)};
			\addplot[line width=1pt,solid,color=blue,mark=square*,dashed,mark options=solid]%
			coordinates {(10, -3.8448) (20, -3.8326) (30, -3.8214) (40, -3.8068) (50, -3.7713) (60, -3.6062) (70, -3.4948) (80, -3.3984) (90, -3.2942)};
			\addplot[line width=1pt,solid,color=red,mark=*,dotted,mark options=solid]%
			coordinates {(10, -3.5156) (20, -3.4955) (30, -3.4807) (40, -3.4685) (50, -3.4567) (60, -3.4449) (70, -3.432) (80, -3.4168) (90, -3.3959)};
		\end{axis}
	\end{tikzpicture}
	\begin{tikzpicture}[scale=0.5]
		\centering
		\begin{axis}[
			legend style={nodes={scale=0.5, transform shape}}, 
			legend image post style={scale=0.5},
			scaled ticks=false, 
			tick label style={/pgf/number format/fixed},
			xlabel={Quantile},
			ylabel={$\lambda_{L}$},
			xmin = 10, xmax = 90,
			ymin = -4.2, ymax = -3.2,
			legend pos=north west,
			]	
			\node at (5mm,52.5mm) {\Large\bfseries\sffamily\selectfont e};
			\addplot[line width=1pt,solid,color=brown,mark=triangle*,solid,mark options=solid]%
			coordinates {(10, -4.0595) (20, -3.9662) (30, -3.8897) (40, -3.8194) (50, -3.7413) (60, -3.6472) (70, -3.5238) (80, -3.4019) (90, -3.2751)};
			\addplot[line width=1pt,solid,color=blue,mark=square*,dashed,mark options=solid]%
			coordinates {(10, -3.8448) (20, -3.8326) (30, -3.8214) (40, -3.8068) (50, -3.7713) (60, -3.6062) (70, -3.4948) (80, -3.3984) (90, -3.2942)};
			\addplot[line width=1pt,solid,color=red,mark=*,dotted,mark options=solid]%
			coordinates {(10, -3.8357) (20, -3.8231) (30, -3.8128) (40, -3.803) (50, -3.792) (60, -3.7789) (70, -3.7558) (80, -3.5958) (90, -3.4108)};
		\end{axis}
	\end{tikzpicture}
	\begin{tikzpicture}[scale=0.5]
		\centering
		\begin{axis}[
			legend style={nodes={scale=0.5, transform shape}}, 
			legend image post style={scale=0.5},
			scaled ticks=false, 
			tick label style={/pgf/number format/fixed},
			xlabel={Quantile},
			ylabel={$\lambda_{L}$},
			xmin = 10, xmax = 90,
			ymin = -4.2, ymax = -3.2,
			legend pos=north west,
			]	
			\node at (5mm,52.5mm) {\Large\bfseries\sffamily\selectfont f};
			\addplot[line width=1pt,solid,color=brown,mark=triangle*,solid,mark options=solid]%
			coordinates {(10, -4.0595) (20, -3.9662) (30, -3.8897) (40, -3.8194) (50, -3.7413) (60, -3.6472) (70, -3.5238) (80, -3.4019) (90, -3.2751)};
			\addplot[line width=1pt,solid,color=blue,mark=square*,dashed,mark options=solid]%
			coordinates {(10, -3.8448) (20, -3.8326) (30, -3.8214) (40, -3.8068) (50, -3.7713) (60, -3.6062) (70, -3.4948) (80, -3.3984) (90, -3.2942)};
			\addplot[line width=1pt,solid,color=red,mark=*,dotted,mark options=solid]%
			coordinates {(10, -3.8198) (20, -3.7966) (30, -3.7777) (40, -3.755) (50, -3.7219) (60, -3.5966) (70, -3.4535) (80, -3.3458) (90, -3.2329)};
		\end{axis}
	\end{tikzpicture}
	\begin{tikzpicture}[scale=0.5]
		\centering
		\begin{axis}[
			legend style={nodes={scale=0.5, transform shape}}, 
			legend image post style={scale=0.5},
			scaled ticks=false, 
			tick label style={/pgf/number format/fixed},
			xlabel={Quantile},
			ylabel={$\lambda_{L}$},
			xmin = 10, xmax = 90,
			ymin = -4.2, ymax = -3.2,
			legend pos=north west,
			]	
			\node at (5mm,52.5mm) {\Large\bfseries\sffamily\selectfont g};
			\addplot[line width=1pt,solid,color=brown,mark=triangle*,solid,mark options=solid]%
			coordinates {(10, -4.1321) (20, -4.0263) (30, -3.945) (40, -3.8807) (50, -3.8095) (60, -3.7254) (70, -3.6252) (80, -3.5276) (90, -3.4168)};
			\addplot[line width=1pt,solid,color=blue,mark=square*,dashed,mark options=solid]%
			coordinates {(10, -3.8703) (20, -3.869) (30, -3.8678) (40, -3.8661) (50, -3.7784) (60, -3.631) (70, -3.5425) (80, -3.466) (90, -3.3722)};
			\addplot[line width=1pt,solid,color=red,mark=*,dotted,mark options=solid]%
			coordinates {(10, -3.4685) (20, -3.4639) (30, -3.4609) (40, -3.4583) (50, -3.456) (60, -3.4534) (70, -3.4508) (80, -3.4477) (90, -3.4436)};
		\end{axis}
	\end{tikzpicture}
	\begin{tikzpicture}[scale=0.5]
		\centering
		\begin{axis}[
			legend style={nodes={scale=0.5, transform shape}}, 
			legend image post style={scale=0.5},
			scaled ticks=false, 
			tick label style={/pgf/number format/fixed},
			xlabel={Quantile},
			ylabel={$\lambda_{L}$},
			xmin = 10, xmax = 90,
			ymin = -4.2, ymax = -3.2,
			legend pos=north west,
			]	
			\node at (5mm,52.5mm) {\Large\bfseries\sffamily\selectfont h};
			\addplot[line width=1pt,solid,color=brown,mark=triangle*,solid,mark options=solid]%
			coordinates {(10, -4.1321) (20, -4.0263) (30, -3.945) (40, -3.8807) (50, -3.8095) (60, -3.7254) (70, -3.6252) (80, -3.5276) (90, -3.4168)};
			\addplot[line width=1pt,solid,color=blue,mark=square*,dashed,mark options=solid]%
			coordinates {(10, -3.8703) (20, -3.869) (30, -3.8678) (40, -3.8661) (50, -3.7784) (60, -3.631) (70, -3.5425) (80, -3.466) (90, -3.3722)};
			\addplot[line width=1pt,solid,color=red,mark=*,dotted,mark options=solid]%
			coordinates {(10, -3.8267) (20, -3.8243) (30, -3.8225) (40, -3.8209) (50, -3.8193) (60, -3.8177) (70, -3.8157) (80, -3.8122) (90, -3.6047)};
		\end{axis}
	\end{tikzpicture}
	\begin{tikzpicture}[scale=0.5]
		\centering
		\begin{axis}[
			legend style={nodes={scale=0.5, transform shape}}, 
			legend image post style={scale=0.5},
			scaled ticks=false, 
			tick label style={/pgf/number format/fixed},
			xlabel={Quantile},
			ylabel={$\lambda_{L}$},
			xmin = 10, xmax = 90,
			ymin = -4.2, ymax = -3.2,
			legend pos=north west,
			]	
			\node at (5mm,52.5mm) {\Large\bfseries\sffamily\selectfont i};
			\addplot[line width=1pt,solid,color=brown,mark=triangle*,solid,mark options=solid]%
			coordinates {(10, -4.1321) (20, -4.0263) (30, -3.945) (40, -3.8807) (50, -3.8095) (60, -3.7254) (70, -3.6252) (80, -3.5276) (90, -3.4168)};
			\addplot[line width=1pt,solid,color=blue,mark=square*,dashed,mark options=solid]%
			coordinates {(10, -3.8703) (20, -3.869) (30, -3.8678) (40, -3.8661) (50, -3.7784) (60, -3.631) (70, -3.5425) (80, -3.466) (90, -3.3722)};
			\addplot[line width=1pt,solid,color=red,mark=*,dotted,mark options=solid]%
			coordinates {(10, -3.8292) (20, -3.8196) (30, -3.8121) (40, -3.8049) (50, -3.7948) (60, -3.7009) (70, -3.5686) (80, -3.4762) (90, -3.3789)};
		\end{axis}
	\end{tikzpicture}
	\caption{$\lambda_{L}$ versus quantile plot for perturbations about $x$ with $\sigma = 1, 0.1$ and $0.01$ from left to right for the $\beta_{\neq}$ case which is given by the solid gold lines with a triangle for each data point, the $\beta_{=}$ case which is given by the dashed blue lines with a square for each data point, and the $\beta_{u}$ case which is given by the dotted red lines with a circle for each data point. Each method and corresponding uncertainty is tested using $\lambda_{L}(\{x^{j} + \Delta x^{j}\})$ with $\sigma = 1$, $\sigma = 0.1$ and $\sigma = 0.01$ from top to bottom for each plot respectively with each $\epsilon$ is picked from a Gaussian distribution.}
	\label{fig:28}
\end{figure}
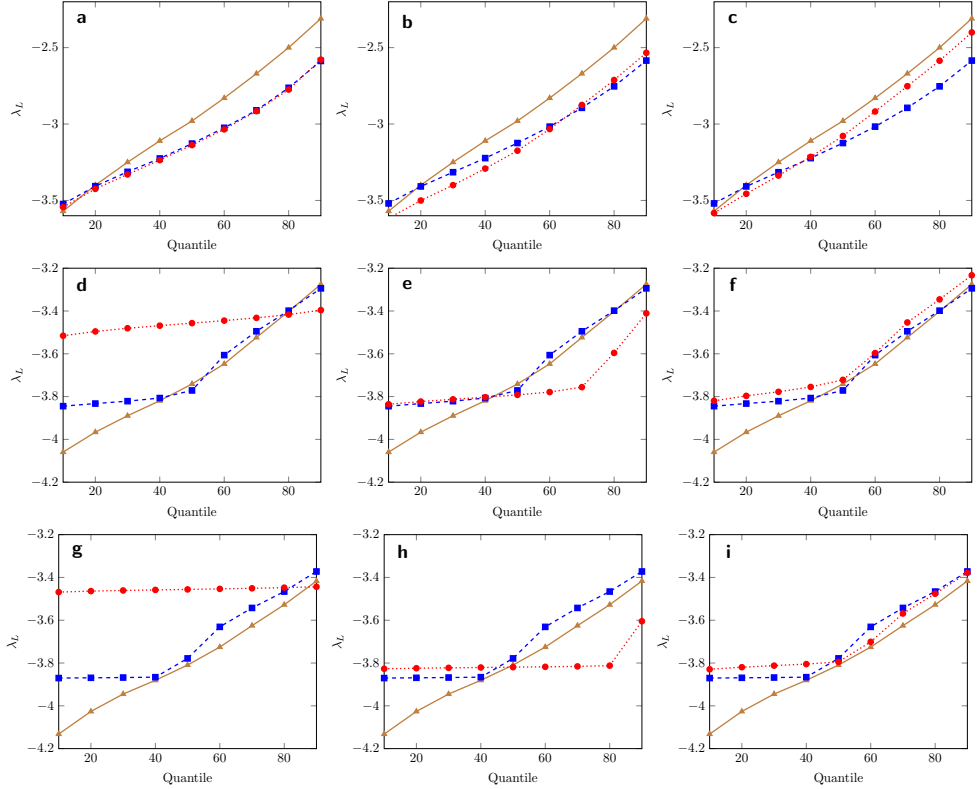

In this section, we analyse each distribution for a given uncertainty to determine which configuration is best and under what circumstances. We look at all aforementioned cases when optimised with various levels of uncertainty and test each method with different degrees of noise. We look at \cref{fig:28} \textbf{a}-\textbf{c} which tests each method with $\lambda_{L}(x^{j} + \Delta x^{j})$ having $\sigma = 1$ when $\{\beta_{i}\}$ is optimised with uncertainty $\sigma = 1$, $\sigma = 0.1$ and $\sigma = 0.01$ respectively. Each plot follows the same trend which is that $\beta_{\neq}$ performs the worst of all methods. This is clear from \cref{fig:28} \textbf{a}-\textbf{c} as the $\beta_{\neq}$ curve is above the $\beta_{=}$ and $\beta_{u}$ curves which mean that for a given quantile, the Lyapunov exponent is greater for the $\beta_{\neq}$ method. This is not surprising given that this method incorporates no uncertainty in the optimisation process and the method finds sharp cusps in the solution space which are not robust to perturbations. The best performing method is the $\beta_{u}$ method for $\sigma = 1$ while $\beta_{=}$ performs better for $\sigma = 0.01$ which suggests that the higher levels of uncertainty in the optimisation process require the $\beta_{u}$ method. However, if there is low uncertainty in the optimisation process then it is sufficient to use the simpler $\beta_{=}$ method even when testing with $\sigma = 1$.

We then look at \cref{fig:28} \textbf{d}-\textbf{f} and \textbf{g}-\textbf{i} which tests each method with $\lambda_{L}(x^{j} + \Delta x^{j})$ having $\sigma = 0.1$ and $\sigma = 0.01$ respectively when $\{\beta_{i}\}$ is optimised with uncertainty $\sigma = 1$, $\sigma = 0.1$ and $\sigma = 0.01$. Although we test these methods with different levels of uncertainty, for both tests they follow a similar trend which is that $\beta_{\neq}$ performs the best initially for low quantiles in all cases. This is to be expected given the relatively low levels of uncertainty when testing which means that there is minimal deviation in each $\beta_{i}$ from its mean in the testing process. At high quantiles, typically above $60$, the $\beta_{u}$ then performs better than the $\beta_{\neq}$ case when optimising with $\sigma = 0.1$. This is a result of the fact that this method is optimised to be robust to perturbations in the power system parameters while also minimising computational complexity to allow an appropriate minimum to be found. However, it is clear that for small quantiles below $40$ when testing the system with $\sigma = 0.1$ or $\sigma = 0.01$, $\beta_{\neq}$ is the best performing method. This demonstrates the advantages of using the $\beta_{\neq}$ method for this particular system state.

\begin{figure}[h]
	\centering
	\begin{tikzpicture}[scale=0.5]
		\centering
		\begin{axis}[
			legend style={nodes={scale=0.5, transform shape}}, 
			legend image post style={scale=0.5},
			scaled ticks=false, 
			tick label style={/pgf/number format/fixed},
			xlabel={Quantile},
			ylabel={$\lambda_{L}$},
			xmin = 10, xmax = 90,
			ymin = -3.6, ymax = -2.2,
			legend pos=north west,
			]
			\node at (5mm,52.5mm) {\Large\bfseries\sffamily\selectfont a};	
			\addplot[line width=1pt,solid,color=brown,mark=triangle*,solid,mark options=solid]%
			coordinates {(10, -3.57) (20, -3.4) (30, -3.25) (40, -3.11) (50, -2.98) (60, -2.83) (70, -2.67) (80, -2.5) (90, -2.31)};
			\addplot[line width=1pt,solid,color=blue,mark=square*,dashed,mark options=solid]%
			coordinates {(10, -3.5235) (20, -3.4068) (30, -3.3127) (40, -3.2254) (50, -3.1278) (60, -3.0247) (70, -2.9105) (80, -2.7628) (90, -2.5872)};
			\addplot[line width=1pt,solid,color=red,mark=*,dotted,mark options=solid]%
			coordinates {(10, -3.5436) (20, -3.4243) (30, -3.3299) (40, -3.2361) (50, -3.1378) (60, -3.0349) (70, -2.9157) (80, -2.7743) (90, -2.5787)};
		\end{axis}
	\end{tikzpicture}
	\begin{tikzpicture}[scale=0.5]
		\centering
		\begin{axis}[
			legend style={nodes={scale=0.5, transform shape}}, 
			legend image post style={scale=0.5},
			scaled ticks=false, 
			tick label style={/pgf/number format/fixed},
			xlabel={Quantile},
			ylabel={$\lambda_{L}$},
			xmin = 10, xmax = 90,
			ymin = -3.6, ymax = -2.2,
			legend pos=north west,
			]	
			\node at (5mm,52.5mm) {\Large\bfseries\sffamily\selectfont b};	
			\addplot[line width=1pt,solid,color=brown,mark=triangle*,solid,mark options=solid]%
			coordinates {(10, -3.57) (20, -3.4) (30, -3.25) (40, -3.11) (50, -2.98) (60, -2.83) (70, -2.67) (80, -2.5) (90, -2.31)};
			\addplot[line width=1pt,solid,color=blue,mark=square*,dashed,mark options=solid]%
			coordinates {(10, -3.5189) (20, -3.4067) (30, -3.312) (40, -3.2238) (50, -3.1251) (60, -3.0195) (70, -2.8961) (80, -2.7533) (90, -2.5866)};
			\addplot[line width=1pt,solid,color=red,mark=*,dotted,mark options=solid]%
			coordinates {(10, -3.5455) (20, -3.428) (30, -3.3249) (40, -3.2356) (50, -3.1422) (60, -3.0458) (70, -2.9218) (80, -2.7787) (90, -2.5889)};
		\end{axis}
	\end{tikzpicture}
	\begin{tikzpicture}[scale=0.5]
		\centering
		\begin{axis}[
			legend style={nodes={scale=0.5, transform shape}}, 
			legend image post style={scale=0.5},
			scaled ticks=false, 
			tick label style={/pgf/number format/fixed},
			xlabel={Quantile},
			ylabel={$\lambda_{L}$},
			xmin = 10, xmax = 90,
			ymin = -4.2, ymax = -3,
			legend pos=north west,
			]	
			\node at (5mm,52.5mm) {\Large\bfseries\sffamily\selectfont c};	
			\addplot[line width=1pt,solid,color=brown,mark=triangle*,solid,mark options=solid]%
			coordinates {(10, -4.1321) (20, -4.0263) (30, -3.945) (40, -3.8807) (50, -3.8095) (60, -3.7254) (70, -3.6252) (80, -3.5276) (90, -3.4168)};
			\addplot[line width=1pt,solid,color=blue,mark=square*,dashed,mark options=solid]%
			coordinates {(10, -3.8726) (20, -3.8726) (30, -3.8726) (40, -3.8726) (50, -3.7822) (60, -3.634) (70, -3.5465) (80, -3.4657) (90, -3.3724)};
			\addplot[line width=1pt,solid,color=red,mark=*,dotted,mark options=solid]%
			coordinates {(10, -3.4665) (20, -3.4627) (30, -3.4601) (40, -3.4577) (50, -3.4556) (60, -3.4534) (70, -3.4512) (80, -3.4487) (90, -3.445)};
		\end{axis}
	\end{tikzpicture}
	\caption{$\lambda_{L}$ versus quantile plot for perturbations about $x$ with $\sigma = 1$ for the $\beta_{=}$, $\beta_{\neq}$ and the $\beta_{u}$ case using the $\lambda_{L}(\{x^{j} + \Delta x^{j}\})$, $\lambda_{L}(\{\beta^{j} + \Delta \beta^{j}\})$ and $\lambda_{L}(\{y^{j} + \Delta y^{j}\})$ tests with $\sigma = 1$ for each plot left to right respectively and each $\epsilon$ is picked from a Gaussian distribution. The solid gold line with triangles for data points corresponds to the $\beta_{\neq}$ case, the dashed blue line with squares for data points corresponds to the $\beta_{=}$ case and the dotted red line with circles for data points corresponds to the $\beta_{u}$ case.}
	\label{fig:27}
\end{figure}
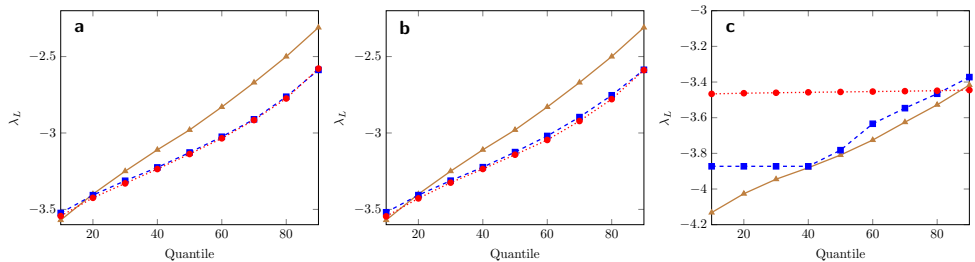

Finally, we also analyse a case where there is either just noise in the damping parameters or power flow parameters which are given in the last two plots in \cref{fig:27}. These cases give similar results to what we saw in previous sections whereby the $\beta_{u}$ method performs the best when tested with $\lambda_{L}(\{x^{j} + \Delta x^{j}\})$ and $\lambda_{L}(\{\beta^{j} + \Delta \beta^{j}\})$ but the worst when tested with $\lambda_{L}(\{y^{j} + \Delta y^{j}\})$. For the $\lambda_{L}(\{y^{j} + \Delta y^{j}\})$ test, the $\beta_{\neq}$ method performs better at low and high quantile values which illustrates that it is not just the level of uncertainty in the system that determines which method is most appropriate but also the parameters that have an associated uncertainty.

We now analyse a worst reasonable case, we define this as the critical $\lambda_{L}$ value denoted as $\lambda_{c}$ that a transmission system operator aims to have a given percentage of $\lambda_{L}$ values fall below. This is essentially the extreme values or tail of the distribution which is of particular interest to a transmission system operator. The $\lambda_{c}$ value varies depending on the method but also on the quantile that the transmission system operator chooses for $\lambda_{c}$. We assume that the transmission system operator would choose a conservative estimate given the negative consequences associated with power grid failures and present $\lambda_{c}$ values for a range of percentages between $95\%-99\%$ for each method.

\begin{table}[h]
	\footnotesize
	\caption{The $\lambda_{c}$ values for the $\beta_{=}$, $\beta_{\neq}$ and $\beta_{u}$ case for a range of worst reasonable case percentages for $\sigma = 1$, $\sigma = 0.1$ and $\sigma = 0.01$.}
	\begin{center}
		\begin{tabular}{|c|c|c|c|c|c|}
			\hline
			\textbf{Method} & \textbf{\textit{$95\%$}} & \textbf{\textit{$96\%$}} & \textbf{\textit{$97\%$}} & \textbf{\textit{$98\%$}} & \textbf{\textit{$99\%$}} \\
			\hline
			\multicolumn{6}{|c|}{$\sigma = 1$}\\
			\hline
			$\beta_{\neq}$ & $-2.17$ & $-2.14$ & $-2.10$ & $-2.06$ & $-2.00$ \\
			\hline
			$\beta_{=}$ & $-2.45$ & $-2.42$ & $-2.38$ & $-2.34$ & $-2.27$ \\
			\hline
			$\beta_{u}$ & $-2.54$ & $-2.50$ & $-2.45$ & $-2.40$ & $-2.31$ \\
			\hline
			\multicolumn{6}{|c|}{$\sigma = 0.1$}\\
			\hline
			$\beta_{\neq}$ & $-3.19$ & $-3.17$ & $-3.14$ & $-3.11$ & $-3.05$ \\
			\hline
			$\beta_{=}$ & $-3.23$ & $-3.21$ & $-3.19$ & $-3.16$ & $-3.11$ \\
			\hline
			$\beta_{u}$ & $-3.57$ & $-3.56$ & $-3.56$ & $-3.55$ & $-3.54$ \\
			\hline
			\multicolumn{6}{|c|}{$\sigma = 0.01$}\\
			\hline
			$\beta_{\neq}$ & $-3.34$ & $-3.32$ & $-3.30$ & $-3.27$ & $-3.22$ \\
			\hline
			$\beta_{=}$ & $-3.31$ & $-3.29$ & $-3.27$ & $-3.25$ & $-3.21$ \\
			\hline
			$\beta_{u}$ & $-3.62$ & $-3.62$ & $-3.62$ & $-3.62$ & $-3.62$ \\
			\hline
		\end{tabular}
		\label{tab7}
	\end{center}
\end{table}

As is clear from \cref{tab7}, for $\sigma = 1$, the methods which generate the best $\lambda_{c}$ are in the following order; $\beta_{u}$, $\beta_{=}$ and $\beta_{\neq}$. This order does not change as we increase the percentage of values $\lambda_{c}$ must be less than. There is a large difference between $\lambda_{c}$ of $\beta_{\neq}$ and the other two methods. This is to be expected given the large levels of uncertainty that when testing $\lambda_{L}$, probe an area beyond the small ridge that the optimal $\beta_{\neq}$ finds in the solution space. There is little difference between $\beta_{=}$ and the best-performing method, although $\beta_{=}$ is significantly easier to compute. As this method is computationally inexpensive and reasonably robust, an argument could be made that this should be the standard method if there are high levels of uncertainty in the system. However, if the most optimal configuration is needed, this justifies using a more complicated, computationally intensive method.

The results in \cref{tab7} for the $\sigma = 0.1$ case are different to that of the higher uncertainty case of $\sigma = 1$ whereby $\beta_{\neq}$ now performs almost as well as $\beta_{=}$. The results for the $\sigma = 0.01$ case have a different ordering of best-performing methods with $\beta_{\neq}$ outperforming $\beta_{=}$ although there is little difference between both methods. The $\beta_{u}$ method is the best-performing method again emphasising the importance of optimising power grid parameters while incorporating uncertainty.

\section{Discussion and Conclusions}
\label{Discussion and Conclusions}
In this paper, we analysed the model of power systems proposed by Molnar et al. and considered how the stability of the system is affected by uncertainty. In \cref{Optimisation Under Uncertainty}, we found that different levels of uncertainty in $\{\beta_{i}\}$ lead to different distributions of $\lambda_{L}$ values. This motivated us to develop an optimisation method that optimises the effective damping parameters of the system and explicitly takes uncertainty into account. Specifically, we used a simulated annealing under uncertainty algorithm where the objective function at each step is based on the average of some $\lambda_{L}$ values associated with a sample of parameter values. In \cref{Statistical Analysis}, we used a quantile-based metric to compare the methods developed by Molnar et al. to the optimisation algorithm we presented. We found that each method had advantages and disadvantages for different levels of uncertainty depending on the quantile of interest. However, the optimisation algorithm that incorporates uncertainty outperformed all other methods at the tail of the distribution for quantiles $95$-$99$.

The level of certainty that the transmission system operator can provide informs the method that should be used for optimising $\{\beta_{i}\}$ parameters. For example, if the uncertainty is unknown for a system parameter, it may be safer to use a conservative estimate by looking at the tail of a $\lambda_{L}$ distribution and choosing a method that is particularly robust to large perturbations even if large perturbations are unlikely to occur. On the other hand, if all system parameters are known with high precision, this opens up the possibility of using the $\beta_{\neq}$ method which would give improved stability with a degree of certainty determined by equipment accuracy. Highly accurate equipment that minimises noise and takes precise measurements of power system parameters to determine their distributions is key to achieving strong stability when using the $\beta_{\neq}$ method. Although the $\beta_{u}$ method gave the best results for high quantile values, it was significantly more computationally expensive than the other methods. Therefore, it is important to understand when it is essential to use this method and when other more computationally simple methods are appropriate.

We used a quantile-based metric to determine which optimisation algorithms gave better stability and under what conditions. This metric is useful because it offers flexibility for transmission system operators when deciding what value a percentage of Lyapunov exponent values must be below. We used this metric as opposed to the mean, mode, median or other standard summary statistics because the tail of the distribution is often of the most importance in power grid stability. Transmission system operators want the Lyapunov exponent to rarely fall beyond their specified bounds which makes the ability to interrogate the tail of the $\lambda_{L}$ distributions appropriate for this application. However, there are other metrics that may be more appropriate than looking at quantiles and an area of future work could be to investigate various other metrics and summary statistics.

An area of future work that may be undertaken is experimental research that would look at determining the distributions of parameters that occur in power grid modelling and quantifying the level of noise and uncertainty. So far, we have assumed white noise but there are various sources of uncertainty and at least some of these sources may have an uncertainty distribution that is not Gaussian. The uncertainties may be different depending on the system of interest but it would be useful to have a consistent benchmark of different levels of uncertainty found in a typical system.

Another area of future work would be to expand on the current model for optimising under uncertainty to make the algorithm more efficient. The current algorithm works well but is computationally expensive given that it is necessary to sample the distributions of $\{\beta_{i}\}$ many times to determine $\bar{\lambda}_{L}^{j}$ when optimising the mean values for a set of $\{\beta_{i}\}$. This may include making use of advanced sampling techniques that do not require the same large number of samples. Alternatively, dimensionality reduction may be a useful tool by reducing the space needed to search for the optimum. The algorithm could also be expanded to use a quantile of $\lambda_{L}$ rather than the mean $\lambda_{L}$ for the objective function. This would have the benefit of incorporating information about the tail of the distribution into the algorithm but would come with an added computational cost as more samples would be needed to ensure an accurate $\lambda_{L}$ is calculated.

In this work, the power flow state is calculated using a steady state instance of the power flow system and then incorporating uncertainty into this instance of the power flow system however, this is a system that changes over time. Parameters from the power flow state that change over time could be incorporated into this model by modelling their behaviour over time and under a variety of conditions. For example, the variation of power generation and demand over a day, week, month and year period could be recorded to provide better distributions of these parameters and a variety of test cases to investigate. Other parameters such as line admittance depend on temperature and under periods of intense stress may heat up and thus these parameters would change. Therefore, an important extension to this model would be to combine the effects of uncertainty over time with that of measurements and noise to create models that are robust to all aspects of uncertainty that may occur in a power system.

Uncertainty plays a significant role in determining the strategy of optimising effective damping parameters to ensure the stability of the grid. It is essential to account for the uncertainty that affects every parameter in the system to model the stability of the grid correctly. Incorporating uncertainty into the optimisation algorithm of the effective damping parameters provides more robust results than previously used methods. The results from the simulated annealing under uncertainty algorithm are optimal at large quantile values which are key to ensuring the stability of the power grid.

\bibliographystyle{siamplain}
\bibliography{references}
\end{document}